# Multi Methods of Matrix Analysis Use for Control and Optimization system in Control Engineering

Moeurn Si Kheang
Beijing Institute of Technology, School of Automation, Beijing, China
E-mail: sikheang2498@163.com

**Abstract**

Matrix analysis plays a crucial role in the field of control engineering, providing a powerful mathematical framework for the analysis and design of control systems. This research report explores various applications of matrix analysis in control engineering, focusing on its contributions to system modeling, stability analysis, controllablity, observability, and optimization.

The report also discusses specific examples and case studies to illustrate the practical significance of matrix analysis in addressing real-world control engineering challenges Analyze controllability. Informally, a system is controllable if we can construct a set of inputs that will drive the system to any given state. Analyze observability. Informally, observability means that by controlling the inputs and watching the outputs of a system we can determine what the states were. Optimal Control is a control method that aims to find the optimal control input to achieve the best performance of the system under certain constraints. This performance index can be the system output, energy consumption, time, etc.

## 1. Introduction

The technological revolution of the past century occurred because people learned to control large systems composed of nature, machines, people, and society[1][2]. When they used and improved devices to help them control such systems, they found that they could control bigger and bigger systems. The design of control devices is called "control engineering or Automation Engineering"[3][4]. Early control devices were mechanical, and their design was mainly intuitive. However, recent control devices are algorithms embedded in computers, and their design is very mathematical, such as modern automobiles have a multitude of devices to help the driver control the vehicle, electronic automatic transmission, power steering, and four-wheel antilock brakes, to mention a few[5][6][7]. For fly vehicle an airplane the pilot needs control devices to translate his manual actions into the large forces required to move the wing control surfaces. Some of the devices he uses amplify his strength and others augment his intelligence[8][9].

For System modeling means that the patterns describing the system change with time and the characteristics of the patterns at any time period are interrelated with those of earlier times. A system is a process that converts inputs to outputs[10][11]. A system accepts inputs and, based on the inputs and its present state, creates outputs. A system has no direct control over its inputs[12]. An example of a system used in everyday life is a traffic light. It accepts inputs, such as pedestrians pushing the walk button or cars driving over sensors, and based on its current state, creates outputs that are the colors of the lights in each direction[13][14].

Defining the state of a system is one of the most important, and often most difficult, tasks in system design. The state of the system is the smallest entity that summarizes the past history of the system[15][16]. The state of the system and the sequence of inputs allows computation of the future states of the system[17]. The state of a system contains all the information needed to calculate future responses without reference to the history of inputs and responses[18].

For analyze stability of systems. Without giving a formal definition we can say that in an unstable system the state can have large variations and small inputs may produce very large outputs[19]. A common example of an unstable system is illustrated by someone pointing the microphone of a public address (PA) system at a speaker[20]; aloud high-pitched tone results[21]. Often instabilities are caused by too much gain. So to quiet the PA system, decrease the gain by pointing the microphone away from the speakers. Discrete systems can also be unstable.

Analyze controllability. Informally, a system is controllable if we can construct a set of inputs that will drive the system to any given state. Analyze observability. Informally, observability means that by controlling the inputs and watching the outputs of a system we can determine what the states were[22][23].

Optimal Control is a control method that aims to find the optimal control input to achieve the best performance of the system under certain constraints. This performance index can be the system output, energy consumption, time, etc. In the theory of optimal control, dynamic programming, variational methods, linear programming, quadratic programming and other mathematical methods are commonly used to solve optimal control problems[24][25].

# 2. System Modeling

## 2.1. The concept of dynamic systems

Let u(t), x(t) and y(t) denote the input, state, and output of a system at time period t. Then the system is represented as the block diagram shown in Figure .

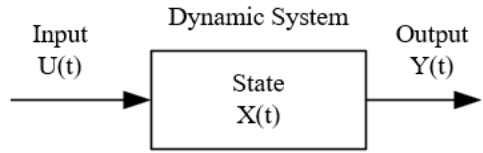

Figure 2.1 Block diagram representation of dynamic systems.

A dynamic system with continuous time scale and state-space description is presented as

$$\dot{X}(t) = f(t, X(t), U(t)) \quad (2.1)$$
$$Y(t) = g(t, X(t)) \quad (2.2)$$

where the first equation is known as the state transition equation and the second relation is known as the output equation. Some authors allow function g to depend also on the input. it is assumed that for all $t \geq 0$,

$$X(t) \in X, \quad U(t) \in U$$

where $X \in \mathbb{R}^n$ and $Y \in \mathbb{R}^m$ are called the state space and input space of the system; furthermore,

$$f : [0, \infty) \times X \times U \to \mathbb{R}^n \text{ and } g : [0, \infty) \times X \to \mathbb{R}^p$$

A dynamic system with discrete time scale is presented as:

$$X(t+1) = f(t, X(t), U(t)) \quad (2.3)$$
$$Y(t) = g(t, X(t)) \quad (2.4)$$

where the first equation is known as the state transition equation and the second relation is known as the output equation. Some authors allow function g to depend also on the input. it is assumed that for all $t = 0,1,2, \ldots$.

$$X(t) \in X, \quad U(t) \in U$$

where $X \in \mathbb{R}^n$ and $Y \in \mathbb{R}^m$ are called the state space and input space of the system; furthermore,

$$f : \mathbb{N} \times X \times U \to \mathbb{R}^n \text{ and } g : \mathbb{N} \times X \to \mathbb{R}^p$$

## 2.2. Continuous linear systems

In this section, dynamic systems of the form:

$$\dot{X}(t) = A(t)X(t) + B(t)U(t), \ X(0) = X_0 \quad (2.5)$$
$$Y(t) = C(t)X(t) \quad (2.6)$$

The differential equation is a special case of the general first-order inhomogeneous linear equation we get,

$$f(t) = B(t)u(t) \quad (2.7)$$

Where matrices A, B and C are constant matrix that show systems state.

**Theorem 2.1**
The general solution of system (2.5) and (2.6) is given by relations

$$X(t) = \phi(t, t_0)X_0 + \int_{t_0}^{t} \phi(t, \tau)B(\tau)U(\tau)d\tau \quad (2.8)$$

And

$$Y(t) = C(t)X(t) =$$
$$C(t)\phi(t, t_0)X_0 + \int_{t_0}^{t} C(t)\phi(t, \tau)B(\tau)U(\tau)d\tau \quad (2.9)$$

These solution formulas are illustrated in the following example.

**Example 2.1**
Now we give the solution of the dynamic system

$$\dot{X}(t) = \begin{pmatrix} 0 & \omega \\ -\omega & 0 \end{pmatrix} X(t) + \begin{pmatrix} 0 \\ 1 \end{pmatrix} U(t), X(0) = \begin{pmatrix} 1 \\ 0 \end{pmatrix}$$
$$Y(t) = (1,1)X(t)$$

Here (1, l) is a row vector.
In example we derived the fundamental matrix:

$$\phi(t, t_0) = \begin{pmatrix} \cos \omega(t - t_0) & \sin \omega(t - t_0) \\ -\sin \omega(t - t_0) & \cos \omega(t - t_0) \end{pmatrix}$$

Therefore, the state variable is

$$X(t) = \begin{pmatrix} \cos \omega t & \sin \omega t \\ -\sin \omega t & \cos \omega t \end{pmatrix} \begin{pmatrix} 1 \\ 0 \end{pmatrix}$$
$$+ \int_{0}^{t} \begin{pmatrix} \cos \omega(t - \tau) \sin \omega(t - \tau) \\ -\sin \omega(t - \tau) \cos \omega(t - \tau) \end{pmatrix} \begin{pmatrix} 0 \\ 1 \end{pmatrix} u(\tau) d\tau$$
$$= \begin{pmatrix} \cos \omega t \\ -\sin \omega t \end{pmatrix} + \int_{0}^{t} \begin{pmatrix} \sin \omega(t - \tau) \\ \cos \omega(t - \tau) \end{pmatrix} u(\tau) d\tau.$$

As a special case, assume that $U(t) = 1$, then the calculations and the state vector become,

$$X(t) = \frac{1}{w} \begin{pmatrix} 1 + (w - 1)\cos \omega t \\ -(w - 1)\sin \omega t \end{pmatrix}$$

The direct input-output relation can be derived as follows:

$$Y(t) = (1,1)X(t)$$
$$= (\cos \omega t - \sin \omega t) + \int_0^t (\sin \omega(t-\tau) + \cos \omega(t-\tau))u(\tau)d\tau$$

In the special case of $U(t) = 1$ we have
$$y(t) = (1,1)\frac{1}{\omega}\begin{pmatrix} 1 + (\omega - 1)\cos \omega t \\ -(\omega - 1)\sin \omega t \end{pmatrix}$$
$$= \frac{1}{\omega}(1 + (\omega - 1)(\cos \omega t - \sin \omega t)).$$

## 2.3. Discrete linear systems

In this section, discrete dynamic systems of the form:
$$X(t+1) = A(t)X(t) + B(t)U(t), \; X(0) = X_0$$
$$Y(t) = C(t)X(t) \quad (2.10)$$

The differential equation is a special case of the general first-order inhomogeneous linear equation we get,
$$f(t) = B(t)u(t) \quad (2.11)$$

Where matrices A, B and C are constant matrix that show systems state.

**Theorem 2.2**
The general solution of system (2.5) and (2.6) is given by relations
$$X(t) = \phi(t,0)X_0 + \sum_{t}^{\tau+1} \phi(t,\tau+1)B(\tau)U(\tau) \quad (2.12)$$

And
$$Y(t) = C(t)X(t) =$$
$$C(t)\phi(t,t_0)X_0 + \sum_t C(t)\phi(t,\tau)B(\tau)U(\tau)d\tau \quad (2.13)$$

These solution formulas are illustrated in the following example.

**Example 2.2**
Now we give the solution of the dynamic system
$$\dot{X}(t) = \begin{pmatrix} 1 & 1 \\ 0 & 1 \end{pmatrix} X(t) + \begin{pmatrix} 0 \\ 1 \end{pmatrix} U(t), X(0) = \begin{pmatrix} 1 \\ 0 \end{pmatrix}$$
$$Y(t) = (1,1)X(t)$$

Here (1, l) is a row vector.
In example we derived the fundamental matrix:
$$\phi(t,t_0) = \begin{pmatrix} \cos \omega(t-t_0) & \sin \omega(t-t_0) \\ -\sin \omega(t-t_0) & \cos \omega(t-t_0) \end{pmatrix}$$

Therefore, the state variable is
$$\phi(t,\tau) = \begin{pmatrix} 1 & t-\tau \\ 0 & 1 \end{pmatrix}$$

As a special case, assume that $U(t) = 1$, then the calculations and the state vector become,

$$X(t) = \begin{pmatrix} 1 & t-0 \\ 0 & 1 \end{pmatrix} \begin{pmatrix} 1 \\ 0 \end{pmatrix} + \sum_{\tau=0}^{t-1} \begin{pmatrix} 1 & t-\tau-1 \\ 0 & 1 \end{pmatrix} \begin{pmatrix} 0 \\ 1 \end{pmatrix} u(\tau)$$
$$= \begin{pmatrix} 1 \\ 0 \end{pmatrix} + \sum_{\tau=0}^{t-1} \begin{pmatrix} t-\tau-1 \\ 1 \end{pmatrix} u(\tau)$$

In the special case of $U(t) = 1$ we have
$$y(t) = (1,1)\mathbf{x}(t) = 1 + \sum_{\tau=0}^{t-1}(t-\tau)u(\tau)$$

In the particular case when $U(t) = 1$, the calculations and the state vector is
$$Y(t) = (1,1)X(t)$$
$$X(t) = \begin{pmatrix} 1 \\ 0 \end{pmatrix} + \sum_{\tau=0}^{t-1} \begin{pmatrix} t-\tau-1 \\ 1 \end{pmatrix} = \begin{pmatrix} (t^2 - t + 2)/2 \\ t \end{pmatrix}$$

## 3. Stability Analysis

Before you begin to format your paper, first write and save the content as a separate text file. Keep your text and graphic files separate until after the text has been formatted and styled. Do not use hard tabs, and limit use of hard returns to only one return at the end of a paragraph. Do not add any kind of pagination anywhere in the paper. Do not number text heads-the template will do that for you.

Finally, complete content and organizational editing before formatting. Please take note of the following items when proofreading spelling and grammar:

### 3.1. Lyapunov stability theory

In this section only time-invariant systems will be considered. Continuous time-invariant systems have the form,
$$\dot{X}(t) = f(X(t)) \quad (3.1)$$

and discrete time-invariant systems are modeled by the difference equation
$$X(t+1) = f(X(t)) \quad (3.2)$$

**Definition 3.1**
(i) An equilibrium point x is stable if there is an $\epsilon > 0$ with the following property: For all $\epsilon_1, 0 < \epsilon_1 < \epsilon_0$, there is an $\epsilon > 0$ such that if $\|\bar{X} - X_0\| < \epsilon$, then $\|\bar{X} - X_0\| < \epsilon_1$ for all $t > t_0$.
(ii) An equilibrium point $\bar{X}$ is asymptotically stable if it is s table and there is an $\epsilon > 0$ such that whenever $\|\bar{X} - X_0\| < \epsilon$, then $X(t) \to \bar{X}(t)$ and $t \to \infty$.
(iii) An equilibrium point $\bar{X}$ is globally asymptotically stable if it is s table and with arbitrary initial state $X_0 \in X, X(t) \to \bar{X}(t)$ as $t \to \infty$.

**Definition 3.2**

A real-valued function $V$ defined on n is called a Lyapunov function, if
(i) $V$ is continuous;
(ii) $V$ has a unique global minimum at $\bar{X}$ with respect to all other points in $\Omega$.
(iii) for any state trajectory $X(t)$ contained in $\Omega$, $V(X(t))$ is nonincreasing in $t$.

**Theorem 3.1** Assume that there exists a Lyapunov function $V$ on the spherical region

$$\Omega = \{X \mid \|\bar{X} - X_0\| < \epsilon\}$$

where $\epsilon > 0$ is given, furthermore $\Omega \subseteq X$ Then the equilibrium is s table.

**Example 3.1**

Consider the differential equation

$$\dot{X}(t) = \begin{pmatrix} 0 & \omega \\ -\omega & 0 \end{pmatrix} X(t) + \begin{pmatrix} 0 \\ 1 \end{pmatrix}$$

we verified that the equilibrium is given as $\bar{X} = (\frac{1}{\omega}, 0)$ and

$$\phi(t, t_0) = \begin{pmatrix} \cos \omega(t - t_0) & \sin \omega(t - t_0) \\ -\sin \omega(t - t_0) & \cos \omega(t - t_0) \end{pmatrix}$$

Assume that the initial state is selected from the neighborhood of the equilibrium, that is,

$$X(0) = \begin{pmatrix} 1/\omega + \alpha \\ \beta \end{pmatrix}$$

where $\alpha$ and $\beta$ in magnitude. The general solution formula implies that

$$X t) = \begin{pmatrix} \cos \omega t \sin \omega t \\ -\sin \omega t \cos \omega t \end{pmatrix} \begin{pmatrix} 1/\omega + \alpha \\ \beta \end{pmatrix}$$

$$+ \int_0^t \begin{pmatrix} \cos \omega(t-\tau) \sin \omega(t-\tau) \\ -\sin \omega(t-\tau) \cos \omega(t-\tau) \end{pmatrix} \begin{pmatrix} 0 \\ 1 \end{pmatrix} d\tau$$

$$= \begin{pmatrix} \left(\frac{1}{\omega} + \alpha\right) \cos \omega t + \beta \sin \omega t \\ -\left(\frac{1}{\omega} + \alpha\right) \sin \omega t + \beta \cos \omega t \end{pmatrix} +$$

$$\begin{pmatrix} \frac{1}{\omega} - \frac{1}{\omega} \cos \omega t \\ \frac{1}{\omega} \sin \omega t \end{pmatrix}$$

The solution implies

$$\mathrm{x}(t) = \begin{pmatrix} \frac{1}{\omega} + \alpha \cos \omega t + \beta \sin \omega t \\ -\alpha \sin \omega t + \beta \cos \omega t \end{pmatrix}$$

and

$$X(t) - \bar{X}(t) = \begin{pmatrix} \alpha \cos \omega t + \beta \sin \omega t \\ -\alpha \sin \omega t + \beta \cos \omega t \end{pmatrix}$$

Then $X(t) - \bar{X}(t)$ solves the homogeneous equation. Simple calculation shows that

$$\|X - \bar{X}\|_2 = \sqrt{\alpha^2 + \beta^2}$$

For the Lyapunov function

$$V(X) = (X - \bar{X})^T (X - \bar{X}) = \|X - \bar{X}\|_2^2$$

Then

$$\frac{d}{dt} V(X(t)) = 2(X - \bar{X})^T \cdot \dot{X} = 2(X - \bar{X})^T (AX + b)$$

$$\frac{d}{dt} V(X) = 2\left(x_1 - \frac{1}{\omega}, x_2\right) \begin{pmatrix} \omega x_2 \\ -\omega x_1 + 1 \end{pmatrix} = 0$$

That is, all conditions of Theorem 3.1 are satisfied, which imply the stability of the equilibrium.
Hence, system is stable processing.

## 3.2. BIBO Stability

BIBO (Bounded Input-Bounded Output)

**Theorem 3.2**

Let $T(t, \tau) = (t_{ij}(t, \tau))$, then the continuous time-variant linear system is BIBO stable if and only if the integral

$$\int_{t_0}^{t} |t_{ij}(t, \tau)| d\tau \qquad (3.4)$$

is bounded for all $t \to to$, $i$ and $j$.

**Theorem 3.3**

Assume that for all eigenvalues $\lambda_i$ of $A$, $RRe\lambda_i < 0$ (or $|\lambda_i| < 1$ ). Then the time-invariant linear continuous (or discrete) system is BIBO stable.

**Example 3.2**

Consider again the continuous system

$$\dot{X}(t) = \begin{pmatrix} 0 & \omega \\ -\omega & 0 \end{pmatrix} X(t) + \begin{pmatrix} 0 \\ 1 \end{pmatrix} U(t)$$

$$Y(t) = (1,1) X(t)$$

In this case, imply that

$$\phi(t, \tau) = \begin{pmatrix} \cos\omega(t-\tau) \sin\omega(t-\tau) \\ -\sin\omega(t-\tau) \cos\omega(t-\tau) \end{pmatrix}$$

Therefore,

$$T(t, \tau) = (1, 1) \begin{pmatrix} \cos\omega(t-\tau) \sin\omega(t-\tau) \\ -\sin\omega(t-\tau) \cos\omega(t-\tau) \end{pmatrix} \begin{pmatrix} 0 \\ 1 \end{pmatrix}$$

$$= \sin\omega(t-\tau) + \cos\omega(t-\tau)$$

and

$$I(t) = \int_0^t |\sin\omega(t-\tau) + \cos\omega(t-\tau)| \, d\tau$$

Hence by selecting $t = \frac{2\pi N}{\omega}$

$$I(t) = \frac{1}{\omega} 4\sqrt{2} N$$

When, $N \to \infty$ the system is BIBO stable.

## 4. Controllability

A dynamic system with initial condition $X(0) = X_0$ is said to be controllable to state $X_1$ at $t_1 (> t_0)$ if there exists an input $U(t)$ such that $X(t_1) = X_1$. This concept is illustrated in Figure 4.1.

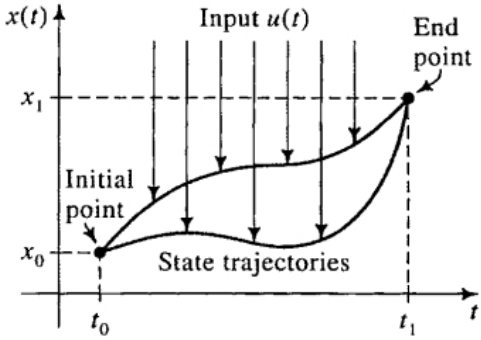

Figure 4.1 Concept of controllability

### 4.1. Continuous Systems

**Theorem 4.1**
The continuous systems is controllable from any initial state $X(t_0) = X_0$ to an arbitrary state $X_1$ at time $t_1 > t_0$ if and only if matrix $W(t_0, t_1)$ is nonsingular.

**Example 4.1**
$$\dot{X}(t) = \begin{pmatrix} 0 & \omega \\ -\omega & 0 \end{pmatrix} X(t) + \begin{pmatrix} 0 \\ 1 \end{pmatrix} U(t), \quad X(0) = \begin{pmatrix} 0 \\ 1 \end{pmatrix}$$

And
$$\phi(t,\tau) = \begin{pmatrix} \cos\omega(t-\tau) & \sin\omega(t-\tau) \\ -\sin\omega(t-\tau) & \cos\omega(t-\tau) \end{pmatrix}$$

Therefore,
$W(0, t_1)$
$$= \int_0^{t_1} \begin{pmatrix} \cos\omega(t-\tau) & \sin\omega(t-\tau) \\ -\sin\omega(t-\tau) & \cos\omega(t-\tau) \end{pmatrix} \begin{pmatrix} 0 \\ 1 \end{pmatrix} (0 \quad 1)$$
$$\begin{pmatrix} \cos\omega(t-\tau) & \sin\omega(t-\tau) \\ -\sin\omega(t-\tau) & \cos\omega(t-\tau) \end{pmatrix} d\tau$$
$$= \begin{pmatrix} \dfrac{t_1}{2} - \dfrac{\sin 2\omega t_1}{4\omega} & \dfrac{\cos 2\omega t_1 - 1}{4\omega} \\ \dfrac{\cos 2\omega t_1 - 1}{4\omega} & \dfrac{t_1}{2} + \dfrac{\sin 2\omega t_1}{4\omega} \end{pmatrix}$$

For all $t_1 > 0$, $W(0, t_1)$ can be written as
$$\dfrac{t_1^2}{4} - \dfrac{\sin^2 2\omega t_1}{16\omega^2} - \dfrac{\cos^2 2\omega t_1 - 2\cos 2\omega t_1 + 1}{16\omega^2}$$

which equals zero if and only if
$$4\omega^2 t_1^2 + 2\cos 2\omega t_1 - 2 = 0$$

Introduce the new variable $\alpha = 2\omega t_1 > 0$ then this equation is equivalent to relation
$$\cos\alpha = 1 - \dfrac{\alpha^2}{2}$$

Consider next function
$$\varphi(\alpha) = \cos\alpha - 1 + \dfrac{\alpha^2}{2}$$

Then easy calculation show that $\varphi(0) = 0$ and all $\alpha > 0$,

$$\varphi'(\alpha) = -\sin\alpha + \alpha > 0$$

Hence, $W(t_0, t_1)$ is usually called the controllability Gramian.

### 4.2. Discrete Systems

**Theorem 4.2**
The discrete system is controllable from initial state $X_0$ to arbitrary state $X_1$ at time $t_1 > 0$ if and only if $W(0, t_1)$ is nonsingular.

**Example 4.2**
Consider the discrete system modeled by difference equation
$$X(t+1) = \begin{pmatrix} 1 & 1 \\ 0 & 1 \end{pmatrix} X(t) + \begin{pmatrix} 0 \\ 1 \end{pmatrix} U(t), \quad X(0) = \begin{pmatrix} 1 \\ 0 \end{pmatrix}$$

And
$$\phi(t,\tau) = \begin{pmatrix} 1 & t-\tau \\ 0 & 1 \end{pmatrix}$$

Therefore,
$W(0, t_1)$
$$= \sum_{\tau=0}^{t_1-1} \begin{pmatrix} 1 & t_1-\tau-1 \\ 0 & 1 \end{pmatrix} \begin{pmatrix} 0 \\ 1 \end{pmatrix} (0 \quad 1) \begin{pmatrix} 1 & 0 \\ t_1-\tau-1 & 1 \end{pmatrix}$$
$$= \begin{pmatrix} \dfrac{t_1(t_1-1)(2t_1-1)}{6} & \dfrac{t_1(t_1-1)}{2} \\ \dfrac{t_1(t_1-1)}{2} & t_1 \end{pmatrix}$$

Determinant $W(0, t_1)$,
$$|W(0, t_1)| = \dfrac{t_1^2}{12}(t_1^2 - 1)$$

Because $W(0, t_1) \neq 0$, hence is usually called the controllability Gramian.

## 5. Observability

**Definition 5.1**
The initial state $X_0$ of a continuous (or discrete) system is said to be observable interval $[t_0, t_1]$ if the trajectories of $U(t)$ and $y(t)$ for $t \in [t_0, t_1]$ (or for $t = 0, 1, 2, \ldots, t_1 - 1$) uniquely determine $X_0$. This concept is illustrated in Figure 5.1.

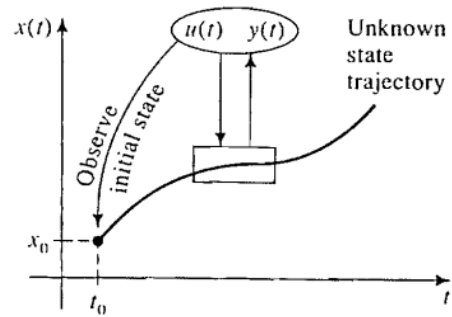

Figure 5.1 Concept of observability

## 5.1. Continuous systems

**Theorem 5.1**

It is possible to determine $X_0$ With in an additive constant vector, which is in $N(M(t_0, t_1))$. If $M(t_0, t_1)$ is nonsingular, then $X_0$ can be determined uniquely.

**Example 5.1**

Consider again the system

$$\dot{X}(t) = \begin{pmatrix} 0 & \omega \\ -\omega & 0 \end{pmatrix} X(t) + \begin{pmatrix} 0 \\ 1 \end{pmatrix} U(t), \quad X(0) = X_0$$

$$Y(t) = (1,1)X(t)$$

And

$$\phi(t,\tau) = \begin{pmatrix} \cos\omega(t-\tau) & \sin\omega(t-\tau) \\ -\sin\omega(t-\tau) & \cos\omega(t-\tau) \end{pmatrix}$$

Therefore,

$$M(0, t_1)$$
$$= \int_0^{t_1} \begin{pmatrix} \cos\omega\tau & -\sin\omega\tau \\ \sin\omega\tau & \cos\omega\tau \end{pmatrix} \begin{pmatrix} 1 \\ 1 \end{pmatrix} (1 \quad 1) \begin{pmatrix} \cos\omega\tau & \sin\omega\tau \\ -\sin\omega\tau & \cos\omega\tau \end{pmatrix} d\tau$$

$$= \begin{pmatrix} t_1 - \dfrac{\sin^2 \omega t_1}{\omega} & \dfrac{\cos\omega t_1 \sin\omega t_1}{\omega} \\ \dfrac{\cos\omega t_1 \sin\omega t_1}{\omega} & t_1 + \dfrac{\sin^2 \omega t_1}{\omega} \end{pmatrix}$$

The non-singularity of this matrix can be examined by computing its determinant, which is

$$t_1^2 - \frac{\sin^4 \omega t_1}{\omega^2} - \left(\frac{\cos\omega t_1 \sin\omega t_1}{\omega}\right)^2 = t_1^2 - \frac{\sin^2 \omega t_1}{\omega^2}$$

Introduce the new variable $\alpha = \omega t_1$, then the determinant equals

$$\frac{\alpha^2}{\omega^2} - \frac{\sin^2 \alpha}{\omega^2} = \frac{1}{\omega^2}(\alpha^2 - \sin^2 \alpha)$$

which is positive for all $\alpha > 0$. Hence, $M(0, t_1)$ is nonsingular for all $t_1 > 0$, therefore, the initial state Xo is observable for all $t_1 > 0$.
Hence the system is completely observable.

## 5.2. Discrete systems

**Theorem 5.2**

It is possible to determine $X_0$ With in an additive constant vector, which is in $N(M(t_0, t_1))$. If $M(t_0, t_1)$ is nonsingular, then $X_0$ can be determined uniquely.

$$M(0, t_1) = \sum_{\tau=0}^{t_1-1} \phi^T(\tau, 0) C^T(\tau) C(\tau) \phi(\tau, 0)$$

**Example 5.1**

Consider again the system

$$\dot{X}(t) = \begin{pmatrix} 0 & \omega \\ -\omega & 0 \end{pmatrix} X(t) + \begin{pmatrix} 0 \\ 1 \end{pmatrix} U(t), \quad X(0) = X_0$$

$$Y(t) = (1,1)X(t)$$

And

$$\phi(t,\tau) = \begin{pmatrix} 1 & t-\tau \\ 0 & 1 \end{pmatrix}$$

Therefore,

$$M(0, t_1) = \sum_{\tau=0}^{t_1-1} \begin{pmatrix} 1 & 0 \\ \tau & 1 \end{pmatrix} \begin{pmatrix} 1 \\ 1 \end{pmatrix} (1 \quad 1) \begin{pmatrix} 1 & \tau \\ 0 & 1 \end{pmatrix}$$

$$= \begin{pmatrix} t_1 & \dfrac{t_1(t_1+1)}{2} \\ \dfrac{t_1(t_1+1)}{2} & \dfrac{t_1(t_1+1)(2t_1+1)}{6} \end{pmatrix}$$

The determinant of $M(0, t_1)$ can be written as

$$|M(0, t_1)| = \frac{t_1^2}{12}(t_1^2 - 1) > 0$$

Hence, for $t_1 \geq 2$, $M(0, t_1)$ is nonsingular and the initial state is observable.

## 6. Optimal Control

### 6.1. Optimal Stability Improvement by State Feedback

Place figures and tables at the top of columns. Avoid placing them in the middle of columns. Large figures

**Theorem 6.1**

Let

$$X(i+1) = AX(i) + BU(i)$$

represent a time-invariant linear discrete-time system. Consider the time invariant control law

$$U(i) = -FX(i)$$

Implies as

$$X(i) = (A - BF)^i X(0)$$

**Example 6.1** (Digital position control system)

The digital positioning system of Example is described

$$X(i+1) = \begin{pmatrix} 1 & 0 \\ 0 & 0.6313 \end{pmatrix} X(i) + \begin{pmatrix} 0.0033 \\ 0.0630 \end{pmatrix} U(i)$$

The system has the characteristic polynomial

$$(Z-1)(Z-0.6313) = Z^2 - 1.6313Z + 0.6313$$

In phase-variable canonical form the system can therefore be represented as

$$Z'(i+1) = \begin{pmatrix} 0 & 1 \\ -0.6313 & 0.6308 \end{pmatrix} X'(i) + \begin{pmatrix} 0 \\ 1 \end{pmatrix} U(i)$$

The transformed state $X'(i)$ is related to the original state $X(i)$ by $X(i) = TX'(i)$, where the matrix $T$ can be found to be

$$T = \begin{pmatrix} 0.00291 & 0.00339 \\ -0.06308 & 0.06308 \end{pmatrix}$$

It is immediately seen that in terms of the transformed state the state dead beat control law is given by

$$U(i) = -(-0.6313, 1.6313)X'(i)$$

In terms of the original state, we have

$$U(i) = -(158.5 \quad 17.33)X(i)$$

**6.2. The Linear Discrete-Time Optimal Regulator Problem**

**Theorem 6.2**

Consider the discrete-time deterministic linear optimal regulator problem. The optimal input is given by,

$$U(i) = -F(i)X(i), \quad i = i_0, i_0 + 1, \ldots, i_1 - 1$$

Where

$$F(i) = \{R(i) + B^T(i)[R_1(i+1) \\ + P(i+1)B(i)\}^{-1}.B^T(i)[R_1(i+1) + P(i+1)]A(i)$$

Here the inverse always exists and

$$R_1(i) = D^T(i)R_s(i)D(i)$$

The sequence of matrices $P(i)$, satisfies the matrix difference equation

$$P(i) = A^T(i)[R_1(i+1) + P(i+1)][A(i) - B(i)F(i)]$$

with the terminal condition

$$P(i) = P_1$$

The value of the criterion achieved with this control law is given by

$$X^T(i_0)P(i_0)X(i_0)$$

**7. Conclusion**

Matrix analysis is a powerful tool in control engineering, playing a crucial role in various areas. One such area is system modeling, where matrices are used to represent the system components and their interactions. The use of matrices helps in simplifying the model and making it more transparent.

System stability is also an important concern in control engineering. Matrices are used to analyze the system's behavior over time, helping engineers determine if the system will remain stable under different conditions. This analysis allows engineers to make necessary adjustments to ensure system stability.

Controllability is the ability to control a system using a specific input. Matrix analysis can be used to assess the controllability of a system by examining the system's eigenvectors and eigenvalues. This analysis provides information on which inputs have the most impact on the system's behavior, making it easier to control.

Observability is the ability to determine the system's state using the output measurements. Matrix analysis is used to assess the observability of a system by studying the system's Jordan canonical form. This analysis helps determine which outputs provide the most information about the system's internal state, improving observability.

Optimal control is concerned with finding the optimal input that will achieve the desired system output with minimum resources. Matrix analysis is used in optimal control by formulating the problem as an optimization problem and solving it using matrix calculus. This approach allows engineers to find the optimal input that will achieve maximum efficiency while adhering to system constraints.

In summary, matrix analysis plays a crucial role in various areas of control engineering, including system modeling, stability analysis, controllability and observability assessments, and optimal control design. These applications have significantly improved the efficiency and performance of control systems.

**References**

[1] Wang Y, Zhang N, Kang C, et al. Standardized matrix modeling of multiple energy systems[J]. IEEE Transactions on Smart Grid, 2017, 10(1): 257-270.

[2] Zhao T, Sun Q H, Li X, et al. A novel transfer matrix-based method for steady-state modeling and analysis of thermal systems[J]. Energy, 2023, 281: 128280.

[3] Luong Q T, Faugeras O D. The fundamental matrix: Theory, algorithms, and stability analysis[J]. International journal of computer vision, 1996, 17(1): 43-75.

[4] Nam H K, Kim Y K, Shim K S, et al. A new eigen-sensitivity theory of augmented matrix and its applications to power system stability analysis[J]. IEEE Transactions on Power Systems, 2000, 15(1): 363-369..

[5] Xu S, Lam J. A survey of linear matrix inequality techniques in stability analysis of delay systems[J]. International Journal of Systems Science, 2008, 39(12): 1095-1113.

[6] Zhang Z, Chen Z, Liu Z. Reachability and controllability analysis of probabilistic finite automata via a novel matrix method[J]. Asian Journal of Control, 2019, 21(6): 2578-2586.

[7] Hahn J, Edgar T F, Marquardt W. Controllability and observability covariance matrices for the analysis and order reduction of stable nonlinear systems[J]. Journal of process control, 2003, 13(2): 115-127.

[8] Klamka J. Controllability of dynamical systems. A survey[J]. Bulletin of the Polish Academy of Sciences. Technical Sciences, 2013, 61(2).

[9] Bei G. Observability analysis for state estimation using Hachtel's augmented matrix method[J]. Electric power systems research, 2007, 77(7): 865-875.

[10] Gou B, Abur A. A direct numerical method for observability analysis[J]. IEEE Transactions on Power Systems, 2000, 15(2): 625-630.

[11] Gou B, Abur A. A direct numerical method for observability analysis[J]. IEEE Transactions on Power Systems, 2000, 15(2): 625-630.

[12] Mouroutsos S G, Sparis P D. Taylor series approach to system identification, analysis and optimal control[J]. Journal of the Franklin Institute, 1985, 319(3): 359-371.

[13] Chen W L, Shih Y P. Analysis and optimal control of time-varying linear systems via Walsh functions[J]. International Journal of Control, 1978, 27(6): 917-932.

[14] Paraskevopoulos P N. Chebyshev series approach to system identification, analysis and optimal control[J]. Journal of the Franklin Institute, 1983, 316(2): 135-157.

[15] Kwakernaak H, Sivan R. Linear optimal control systems[M]. New York: Wiley-interscience, 1972.

[16] Vadali S R, Kim E S .Feedback control of tethered satellites using Lyapunov stability theory[J].Journal of Guidance Control & Dynamics, 2015, 14(4):729-735.DOI:10.2514/3.20706.

[17] Menguc E C, Acir N .Lyapunov Stability Theory Based Adaptive Filter Algorithm for Noisy Measurements[C]//Uksim International Conference on Computer Modelling & Simulation.IEEE, 2013.DOI:10.1109/UKSim.2013.50.

[18] Kim D Y, Park J B, Joo Y H .Reducing conservativeness in stability conditions of affine fuzzy systems using fuzzy Lyapunov function[C]//Fuzzy Systems (FUZZ), 2013 IEEE International Conference on.IEEE, 2013.DOI:10.1109/FUZZ-IEEE.2013.6622411.

[19] Wu X, Wang Y, Huang L ,et al.Robust stability analysis of delayed Takagi-Sugeno fuzzy Hopfield neural networks with discontinuous activation functions[J].Cognitive Neurodynamics, 2010, 4(4):347-354.DOI:10.1007/s11571-010-9123-z.

[20] Russell D L .Controllability and Stabilizability Theory for Linear Partial Differential Equations: Recent Progress and Open Questions[J].Siam Review, 1978, 20(4):639-739.DOI:10.1137/1020095.

[21] Liu Y Y, Slotine J J ,Barabási, Albert-László.Controllability of complex networks[J].Nature, 2011, 473(7346):167.DOI:10.1038/nature10011.

[22] Hermann R, Krener A .Nonlinear controllability and observability[J].IEEE Trans.automat.contr, 1977, 22(5):728-740.DOI:10.1109/TAC.1977.1101601.

[23] Russell D L .A unified boundary controllability theory for hyperbolic and parabolic distributed parameter systems[C]//IEEE Conference on Decision & Control.IEEE, 1972.DOI:10.1109/CDC.1972.269020.

[24] Cheng D, Qi H .Controllability and observability of Boolean control networks[J].Control Theory & Applications, 2009, 45(7):1659-1667.DOI:10.1016/j.automatica.2009.03.006.

[25] Ramadge P J .Observability of discrete event systems[C]//IEEE Conference on Decision & Control.IEEE, 1986.DOI:10.1109/CDC.1986.267551.